\documentclass[11pt,dvipdfm]{article}
\usepackage{amsfonts}
\usepackage{graphicx}
\usepackage{amsfonts, amsmath, amssymb, latexsym, epsfig}
\usepackage{mathrsfs}
\newcommand{\blst}{\begin{trivlist}}
\newcommand{\elst}{\end{trivlist}}

\newtheorem{thm}{Theorem}[section]
\newtheorem{prop}[thm]{Proposition}
\newtheorem{cor}[thm]{Corollary}
\newtheorem{lem}[thm]{Lemma}
\newtheorem{conj}[thm]{Conjecture}
\newtheorem{exa}[thm]{Example}
\newtheorem{defn}[thm]{Definition}

\newcommand{\ben}{\begin{enumerate}}
\newcommand{\een}{\end{enumerate}}
\newcommand{\ble}{\begin{lem}}
\newcommand{\ele}{\end{lem}}
\newcommand{\bth}{\begin{thm}}
\renewcommand{\eth}{\end{thm}}
\newcommand{\bpr}{\begin{prop}}
\newcommand{\epr}{\end{prop}}
\newcommand{\bco}{\begin{cor}}
\newcommand{\eco}{\end{cor}}
\newcommand{\bcon}{\begin{conj}}
\newcommand{\econ}{\end{conj}}
\newcommand{\bde}{\begin{defn}}
\newcommand{\ede}{\end{defn}}
\newcommand{\bex}{\begin{exa}}
\newcommand{\eex}{\end{exa}}
\newcommand{\barr}{\begin{array}}
\newcommand{\earr}{\end{array}}
\newcommand{\btab}{\begin{tabular}}
\newcommand{\etab}{\end{tabular}}
\newcommand{\beq}{\begin{equation}}
\newcommand{\eeq}{\end{equation}}

\newcommand{\bea}{\begin{eqnarray*}}
\newcommand{\eea}{\end{eqnarray*}}
\newcommand{\beaa}{\begin{eqnarray}}
\newcommand{\eeaa}{\end{eqnarray}}
\newcommand{\bce}{\begin{center}}
\newcommand{\ece}{\end{center}}
\newcommand{\bpi}{\begin{picture}}
\newcommand{\epi}{\end{picture}}
\newcommand{\bfi}{\begin{figure} \begin{center}}
\newcommand{\efi}{\end{center} \end{figure}}

\newcommand{\bsl}{\begin{slide}{}}
\newcommand{\esl}{\end{slide}}

{}
{}
{}
{}
{}
{}
{}
{}
{}
{}
{}
{}
{}
{}
{}
{}
{}
{}
{}
{}
{}
{}
{}
{}
{}
{}
{}
{}
\newenvironment{proof}{
\par
\noindent {\bf Proof.}\rm}{\mbox{}\hfill\rule{0.5em}{0.809em}\par}
\begin{document}
\title{Refinements of Dyck Paths with Flaws}

\author{
Jun Ma$^{a,}$\thanks{Email address of the corresponding author:
majun@math.sinica.edu.tw}
 \and Yeong-Nan Yeh $^{b,}$\thanks{Partially supported by NSC 96-2115-M-001-005}}

\date{}
\maketitle \vspace{-1cm} \bce \footnotesize
 $^{a,b}$ Institute of Mathematics, Academia Sinica, Taipei, Taiwan\\
\ece

\thispagestyle{empty}\vspace*{.4cm}

\begin{abstract}
The classical Chung-Feller theorem \cite{CF} tells us that the
number of Dyck paths of length $n$ with $m$ flaws is the $n$-th
Catalan number and independent on $m$. In this paper, we consider
the refinements of Dyck paths with flaws by four parameters, namely
peak, valley, double descent and double ascent. Let ${p}_{n,m,k}$ be
the number of all the Dyck paths of semi-length $n$ with $m$ flaws
and $k$ peaks. First, we derive the reciprocity theorem for the
polynomial $P_{n,m}(x)=\sum\limits_{k=1}^np_{n,m,k}x^k$. Then we
find the Chung-Feller properties for the sum of $p_{n,m,k}$ and
$p_{n,m,n-k}$. Finally, we provide a Chung-Feller type theorem for
Dyck paths of length $n$ with $k$ double ascents: the number of all
the Dyck paths of semi-length $n$ with $m$ flaws and $k$ double
ascents is equal to the number of all the Dyck paths that
 have semi-length $n$, $k$ double ascents and never pass below
the $x$-axis, which is counted by the Narayana number. Let
${v}_{n,m,k}$ (resp. $d_{n,m,k}$) be the number of all the Dyck
paths of semi-length $n$ with $m$ flaws and $k$ valleys (resp.
double descents). Some similar results are derived.
\end{abstract}

\section{Introduction}
Let $\mathbb{Z}$ denote the set of the integers. We consider
$n$-Dyck paths in the plane $\mathbb{Z}\times \mathbb{Z}$ using {\it
up $(1,1)$} and {\it down $(1,-1)$} steps that go from the origin to
the point (2n,0). We say $n$ the \emph{semilength} because there are
$2n$ steps. Define $\mathcal{L}_n$ as the set of all $n$-Dyck paths.
Let $\mathcal{L}=\bigcup\limits_{n\geq 0}\mathcal{L}_n$. A $n$-{\it
flawed} path is a $n$-Dyck path that contains some steps under the
$x$-axis. The number of $n$-Dyck path that never pass below the
$x$-axis is the $n$-th Catalan number
$c_n=\frac{1}{n+1}{2n\choose{n}}$. Such paths are called the {\it
Catalan paths of length $n$}. The generating function
$C(z):=\sum_{n\ge 0}c_n z^n$ satisfies the functional equation
$C(z)=1+zC(z)^2$ and $C(z)=\frac{1-\sqrt{1-4z}}{2z}$ explicitly.

A Dyck path is called a $(n,m)$-{\it flawed} path if it contains $m$
up steps under the $x$-axis and its semilength is $n$. Clearly,
$0\leq m\leq n$. Let $\mathcal{L}_{n,m}$ denote the set of all  the
$(n,m)$-flawed
 paths and $l_{n,m}=|\mathcal{L}_{n,m}|$. The classical Chung-Feller
 theorem \cite{CF} says that $l_{n,m}=c_n$ for $0\leq m\leq n$.

We can consider an $(n,m)$-flawed path $P$ as a word of $2n$ letters
using only $U$ and $D$. In this word, let $P_i$ denote the $i$-th
($1\le i\le 2n$) letter from the left. If a joint node in the Dyck
path is formed by a up step followed by a down step, then this node
is called a {\it peak}; if a joint node in the Dyck path is formed
by a down step followed by a up step, then this node is called a
{\it valley}; if a joint node in the Dyck path is formed by a up
step followed by a up step, then this node is called a {\it double
ascent}; if a joint node in the Dyck path is formed by a down step
followed by a down step, then this node is called a {\it double
descent}.

Define ${\mathcal{P}}_{n,m,k}$ (resp. $\mathcal{{V}}_{n,m,k}$) as
the set of all  the $(n,m)$-flawed paths with $k$ peaks (resp.
valleys). Let $p_{n,m,k}=|\mathcal{P}_{n,m,k}|$ and
$v_{n,m,k}=|\mathcal{V}_{n,m,k}|$. We also define
$\mathcal{A}_{n,m,k}$ (resp. $\mathcal{D}_{n,m,k}$) as the set of
$(n,m)$-flawed path with $k$ double ascents (resp. $k$ double
descents). Let $a_{n,m,k}=|\mathcal{A}_{n,m,k}|$ and
$d_{n,m,k}=|\mathcal{D}_{n,m,k}|$. Let $\varepsilon$ be a mapping
from the set $\{U,D\}$ to itself such that $\varepsilon(U)=D$ and
$\varepsilon(D)=U$. Furthermore, for any path $P=P_1P_2\ldots
P_{2n}\in\mathcal{P}_{n,m,k}$, let
$\phi(P)=\varepsilon(P_1)\varepsilon(P_2)\ldots\varepsilon(P_{2n})$.
It is easy to see that $\phi$ is a bijection between the sets
$\mathcal{P}_{n,m,k}$ and $\mathcal{V}_{n,n-m,k}$. For any
$P=P_1P_2\ldots P_{2n}\in\mathcal{A}_{n,m,k}$, let
$\psi(P)=\varepsilon(P_{2n})\varepsilon(P_{2n-1})\ldots
\varepsilon(P_1)$. Clearly, $\psi$ is a bijection from the set
$\mathcal{A}_{n,m,k}$ to the set $\mathcal{D}_{n,m,k}$. Hence, in
this paper, we focus on the polynomials
$P_{n,m}(x)=\sum\limits_{k=1}^np_{n,m,k}x^k$ and
$A_{n,m}(x)=\sum\limits_{k=0}^{n-1}a_{n,m,k}x^k$. Table $1$ shows
the polynomials $P_{n,m}(x)$ for small values of $n$ and $m$.
$$\begin{array}{||l|l|l|l||} \hline (n,m)&P_{n,m}(x)&(n,m)&P_{n,m}(x)\\
\hline (1,0)&x&(5,0)&x^5+10x^4+20x^3+10x^2+x\\
\hline (1,1)&1&(5,1)&5x^4+20x^3+15x^2+2x\\
\hline (2,0)&x^2+x& (5,2)&4x^4+18x^3+17x^2+3x\\
\hline (2,1)&2x&(5,3)&3x^4+17x^3+18x^2+4x\\
\hline (2,2)&x+1& (5,4)&2x^4+15x^3+20x^2+5x\\
\hline (3,0)&x^3+3x^2+x& (5,5)&x^4+10x^3+20x^2+10x+1\\
\hline (3,1)&3x^2+2x& (6,0)&x^6+15x^5+50x^4+50x^3+15x^2+x\\
\hline (3,2)&2x^2+3x& (6,1)&6x^5+40x^4+60x^3+24x^2+2x\\
\hline (3,3)&x^2+3x+1& (6,2)&5x^5+35x^4+60x^3+29x^2+3x\\
\hline (4,0)&x^4+6x^3+6x^2+x& (6,3)&4x^5+32x^4+60x^3+32x^2+4x\\
\hline (4,1)&4x^3+8x^2+2x& (6,4)&3x^5+29x^4+60x^3+35x^2+5x\\
\hline (4,2)&3x^3+8x^2+3x& (6,5)&2x^5+24x^4+60x^3+40x^2+6x\\
\hline (4,3)&2x^3+8x^2+4x& (6,6)&x^5+15x^4+50x^3+50x^2+15x+1\\
\hline (4,4)&x^3+6x^2+6x+1&&\\
\hline
\end{array}
$$\begin{center} Table $1.$
The polynomials $D_{n,m}(x)$ for small values of $n$ and
$m$.\end{center} From the classical Chung-Feller theorem, we have
$P_{n,m}(1)=A_{n,m}(1)=c_n$ for $0\leq m\leq n$. The classical
Chung-Feller theorem was proved by using analytic method in
\cite{CF}. T.V.Narayana \cite{N} showed the theorem by combinatorial
methods. S.P.Eu et al. \cite{EFY2} studied the theorem by using the
Taylor expansions of generating functions and proved a refinement of
this theorem. Y.M. Chen \cite{C} revisited the theorem by
establishing a bijection. Recently, Shu-Chung Liu et al. \cite{Liu}
use an unify algebra approach to prove chung-Feller theorems for
Dyck path and Motzkin path and develop a new method to find some
combinatorial structures which have  the Chung-Feller property.
However, the macroscopical structures should be supported by some
microcosmic structures. We want to find the Chung-Feller phenomenons
in the more exquisite structures.

Richard Stanley's book \cite{SRP4}, in the context of rational
generating functions, devotes an entire section to exploring the
relationships (called reciprocity relationships) between positively-
and nonpositively-indexed terms of a sequence. First, we give the
reciprocity theorem for the polynomial $P_{n,m}(x)$. Particularly,
we prove that the number of Dyck paths of semi-length $n$ with $m$
flaws and $k$ peaks is equal to the number of Dyck paths of
semi-length $n$ with $n-m$ flaws and $n-k$ peaks.

One observes that the sum of $p_{n,m,k}$ and $p_{n,m,n-k}$ are
independent on $m$ for any $1\leq m\leq n-1$ and $1\leq k\leq
\lfloor\frac{n}{2}\rfloor$ in Table $1$. This is proved in Theorem
\ref{chungfeller} by using the algebra methods. Given $n$ and $k$,
we also show that the polynomials $A_{n,m}(x)$ have the Chung-Feller
property on $m$. Particularly, we conclude that the number of all
the Dyck paths of semi-length $n$ with $m$ flaws and $k$ double
ascents is equal to the number of all the Dyck paths that
 have semi-length $n$, $k$ double ascents and never pass below
the $x$-axis, which is counted by the Narayana number. So, the
Classical Chung-Feller theorem can be viewed as the direct corollary
of this result.

This paper is organized as follows. In Section $2$, we will prove
the reciprocity theorem for the polynomial $P_{n,m}(x)$. In Section
$3$, we will show that $p_{n,m,k}+p_{n,m,n-k}$ have the Chung-Feller
property on $m$ for any $1\leq m\leq n-1$ and $1\leq k\leq
\lfloor\frac{n}{2}\rfloor$. In Section $4$, we will prove that the
polynomials $A_{n,m}(x)$ have the Chung-Feller property on $m$.

\section{The reciprocity theorem for the polynomial $P_{n,m}(x)$ }
In this section, first, define the generating functions
$P_m(x,z)=\sum\limits_{n\geq m}P_{n,m}(x)z^n$. When $m=0$,
$p_{n,0,k}=\frac{1}{k}{n-1\choose{k-1}}{n\choose{k-1}}$ is the
Narayana numbers. It is well known that
$$P_0(x,z)=1+P_0(x,z)z\left[x+P_0(x,z)-1\right],$$ equivalently,
$$P_0(x,z)=\displaystyle{\frac{1+(1-x)z-\sqrt{1-2(1+x)z+(1-x)^2z^2}}{2z}}.$$
Similarly, let $V_m(x,z)=\sum\limits_{n,k\geq 0}v_{n,m,k}x^kz^n$. It
is easy to obtain
$$V_0(x,z)=1+z+z[V_0(x,z)-1](1+xV_0(x,z)),$$
equivalently,
$$V_0(x,z)=\displaystyle{\frac{1-(1-x)z-\sqrt{1-2(1+x)z+(1-x)^2z^2}}{2zx}}.$$
In fact, we have $v_{n,0,k}=p_{n,0,k+1}$ since the number of the
valleys is equal to the number of the peaks minus one for each
Catalan path. So, $P_n(x,z)=V_0(x,z)$.

Now, let $P(x,y,z)=\sum\limits_{n\geq
0}\sum\limits_{m=0}^{n}\sum\limits_{k=1}^{n}p_{n,m,k}x^ky^mz^n.$ Let
$P\in \mathcal{L}$ contain some step over $x$-axis. We decompose $P$
into $P_1UP_2DP_3$, where $U$ and $D$ are the first up and down
steps leaving and returning to $x$-axis and on $x$-axis
respectively. Note that all the steps of $P_1$ are below $x$-axis,
$P_2$ is a Catalan path and $P_3\in \mathcal{D}$. If
$P_2=\emptyset$, then we get a peak $UD$. So, we obtain the
following lemma.
\begin{lem}\label{generating}
$$P(x,y,z)=V_0(x,yz)\left\{1+z[x+P_0(x,z)-1]P(x,y,z)\right\}.$$Equivalently,
\begin{eqnarray*}
P(x,y,z)=\displaystyle{\frac{2}{\sqrt{f(x,z)}+\sqrt{f(x,yz)}+(1-x)(1-y)z}}
\end{eqnarray*}where $f(x,y)=1-2(1+x)y+(1-x)^2y^2$.
\end{lem}

We state the reciprocity relationships for the polynomials
$P_{n,m}(x)$ as the following theorem.

\begin{thm}\label{symetric} Let $n\geq 1$.
$P_{n,m}(x)=x^nP_{n,n-m}(\frac{1}{x})$ for all $0\leq m\leq n$.
Equivalently, $p_{n,m,k}=p_{n,n-m,n-k}$
\end{thm}
\begin{proof} Let $f(x,y)=1-2(1+x)y+(1-x)^2y^2$. Note that (1) $f(x^{-1},xyz)=f(x,yz)$; (2) $f(x^{-1},xz)=f(x,z)$;
and (3) $(1-x^{-1})(1-y^{-1})xyz=(1-x)(1-y)z$.

By Lemma \ref{generating}, we have
$$P(x,y,z)=P(x^{-1},y^{-1},xyz).$$ Since
$P(x,y,z)=1+\sum\limits_{n\geq
1}\sum\limits_{m=0}^nP_{n,m}(x)y^mz^n$, we have
\begin{eqnarray*}P(x^{-1},y^{-1},xyz)&=&1+\sum\limits_{n\geq
1}\sum\limits_{m=0}^nP_{n,m}(\frac{1}{x})y^{-m}(xyz)^n\\
&=&1+\sum\limits_{n\geq
1}\sum\limits_{m=0}^nx^nP_{n,m}(\frac{1}{x})y^{n-m}z^n\\
&=&1+\sum\limits_{n\geq
1}\sum\limits_{m=0}^nx^nP_{n,n-m}(\frac{1}{x})y^{m}z^n.
\end{eqnarray*}
This implies $P_{n,m}(x)=x^nP_{n,n-m}(\frac{1}{x})$ for all $0\leq
m\leq n$. Comparing the coefficients on the sides of the identity,
we derive $p_{n,m,k}=p_{n,n-m,n-k}$.
\end{proof}

Recall that $v_{n,m,k}$ is the number of Dyck paths of semi-length
$n$ with $m$ flaws and $k$ valleys and $v_{n,m,k}=p_{n,n-m,k}$.

\begin{cor} Let $n\geq 1$. Then $v_{n,m,k}=v_{n,n-m,n-k}$.
\end{cor}

\section{The refinement of $(n,m)$-flawed paths obtained by peak}
In this section, we will consider the refinement of $(n,m)$-flawed
paths obtained by peak and prove that the values of
$p_{n,m,k}+p_{n,m,n-k}$ have the Chung-Feller property on $m$ for
any $1\leq m\leq n-1$ and $1\leq k\leq \lfloor\frac{n}{2}\rfloor$.
\begin{lem}\label{flawed=1}
$$P_{1}(x,z)=(1+z-xz)P_0(x,z)-1.$$ Furthermore, we have
$$p_{n,1,k}=\frac{2(n-k)}{n(n-1)}{n\choose{k-1}}{n\choose{k}}$$
for any $n\geq 2$.
\end{lem}
\begin{proof} Let $P$ be a Dyck path containing exact one up step
under the $x$-axis. Then we can decompose the path $P$ into
$P_1DUP_2$, where $P_1$ and $P_2$ are both Catalan paths. So,
$P_1(x,z)=z[P_0(x,z)]^2$. Hence, we have
$P_{1}(x,z)=(1+z-xz)P_0(x,z)-1$ since
$P_0(x,z)=1+P_0(x,z)z\left[x+P_0(x,z)-1\right].$

Note that $P_0(x,z)=1+\sum\limits_{n\geq
1}\sum\limits_{k=1}^{n}p_{n,0,k}x^kz^n$, where
$p_{n,0,k}=\frac{1}{k}{n-1\choose{k-1}}{n\choose{k-1}}$. Therefore,
\begin{eqnarray*}p_{n,1,k}&=&\frac{1}{k}{n-1\choose{k-1}}{n\choose{k-1}}+\frac{1}{k}{n-2\choose{k-1}}{n-1\choose{k-1}}-\frac{1}{k-1}{n-2\choose{k-2}}{n-1\choose{k-2}}\\
&=&\frac{2(n-k)}{n(n-1)}{n\choose{k-1}}{n\choose{k}}.\end{eqnarray*}
\end{proof}
\begin{thm}\label{chungfeller} Let $n$ be an integer with $n\geq 1$ and $1\leq k\leq
\left\lfloor\frac{n}{2}\right\rfloor$. Then
\begin{eqnarray*}p_{n,m,k}+p_{n,m,n-k}=p_{n,m,k}+p_{n,n-m,k}=\frac{2(n+2)}{n(n-1)}{n\choose{k-1}}{n\choose{k+1}}\end{eqnarray*}
for any $1\leq m\leq n-1$.
\end{thm}
\begin{proof} Theorem \ref{symetric} implies that
$p_{n,m,k}+p_{n,m,n-k}=p_{n,m,k}+p_{n,n-m,k}$. We consider the
generating function $R(x,y,z)=\sum\limits_{n\geq
1}\sum\limits_{m=1}^{n-1}\sum\limits_{k=1}^{n-1}(p_{n,m,k}+p_{n,n-m,k})x^ky^mz^n.$
It is easy to see
\begin{eqnarray*}R(x,y,z)&=&P(x,y,z)+P(x,y^{-1},yz)+2\\
&&-[V_0(x,z)+V_0(x,yz)]-[P_0(x,z)+P_0(x,yz)]\\\end{eqnarray*} Let
$\alpha(x,z)=\displaystyle{\frac{1+x-(1-x)z}{x}P_0(x,z)-\frac{P_0(x,z)}{V_0(x,z)}-\frac{1}{x}}$.

 Then
 \begin{eqnarray*}R(x,y,z)=
\displaystyle{\frac{y\alpha(x,z)-\alpha(x,yz)}{1-y}}.
 \end{eqnarray*}

 Suppose $\alpha(x,z)=\sum\limits_{n\geq
 1}\sum\limits_{k=1}^{n-1}a_{k,n}x^kz^n$. Then
\begin{eqnarray*}R(x,y,z)&=&\displaystyle{\sum\limits_{n\geq
1}\sum\limits_{k=1}^{n-1}a_{k,n}x^kz^n\frac{y(1-y^{n-1})}{1-y}}\\
&=&\displaystyle{\sum\limits_{n\geq
1}\sum\limits_{m=1}^{n-1}\sum\limits_{k=1}^{n-1}a_{k,n}x^ky^mz^n}.
\end{eqnarray*}
Hence, given $n\geq 1$ and $1\leq k\leq
\left\lfloor\frac{n}{2}\right\rfloor$, we have
$p_{n,m,k}+p_{n,n-m,k}=p_{n,m,k}+p_{n,m,n-k}=a_{k,n}$ for all $1\leq
m\leq n-1$. By Lemma \ref{flawed=1}, we have
\begin{eqnarray*}
p_{n,m,k}+p_{n,n-m,k}&=&p_{n,1,k}+p_{n,1,n-k}\\
&=&\frac{2(n+2)}{n(n-1)}{n\choose{k-1}}{n\choose{k+1}}.
\end{eqnarray*}
\end{proof}
\begin{cor} Let $n$ be an integer with $n\geq 1$. Then
\begin{eqnarray*}p_{2n,m,n}=\displaystyle{\frac{1}{2n-1}{2n\choose{n-1}}{2n\choose{n}}}\end{eqnarray*}
for any $1\leq m\leq 2n-1$.
\end{cor}

Note that $v_{n,m,k}=p_{n,n-m,k}$. We obtain the following
corollaries.

\begin{cor} Let $n$ be an integer with $n\geq 1$ and $1\leq k\leq
\left\lfloor\frac{n}{2}\right\rfloor$. Then
\begin{eqnarray*}v_{n,m,k}+v_{n,m,n-k}=v_{n,m,k}+v_{n,n-m,k}=\frac{2(n+2)}{n(n-1)}{n\choose{k-1}}{n\choose{k+1}}\end{eqnarray*}
for any $1\leq m\leq n-1$.
\end{cor}

\begin{cor} Let $n$ be an integer with $n\geq 1$. Then
\begin{eqnarray*}v_{2n,m,n}=\displaystyle{\frac{1}{2n-1}{2n\choose{n-1}}{2n\choose{n}}}\end{eqnarray*}
for any $1\leq m\leq 2n-1$.
\end{cor}

In the following theorem, we derive a recurrence relation for the
polynomial $P_{n,m}(x)$.
\begin{thm} For any $m,r\geq 0$, we have
$$
P_{m+r,m}(x)=\left\{\begin{array}{lll} 1&\text{if}&(m,r)=(0,0)\\
\sum\limits_{k=1}^m\frac{1}{k}{m-1\choose{k-1}}{m\choose{k-1}}x^{k-1}&\text{if}&r=0\text{ and }m\geq 1\\
x\sum\limits_{i=0}^m\sum\limits_{j=0}^{r-1}P_{m-i,m-i}(x)P_{r-j-1,r-j-1}(x)P_{j+i,i}(x)&\text{if}&r\geq
1
 \end{array}\right.
$$
\end{thm}
\begin{proof} It is trivial for the case with $r=0$. We only
consider the case with $r\geq 1$. Note that
$x+P_0(x,z)-1=xV_0(x,z)$. Lemma \ref{generating} tells us that
\begin{eqnarray}P(x,y,z)=V_0(x,yz)+xzV_0(x,z)V_0(x,yz)P(x,y,z).\end{eqnarray}
It is well know that $V_0(x,z)=\sum\limits_{n\geq 0}b_n(x)z^n$,
where $b_0(x)=1$ and
$b_n(x)=\sum\limits_{k=1}^n\frac{1}{k}{n-1\choose{k-1}}{n\choose{k-1}}x^{k-1}$
for all $n\geq 1$. Comparing the coefficients of $y^m$ on both side
of Identity (1), we get
\begin{eqnarray}P_m(x,z)=b_m(x)z^m+xzV_0(x,z)\sum\limits_{i=0}^mP_{i}(x,z)b_{m-i}(x)z^{m-i}.\end{eqnarray}
Finally, since $P_m(x,z)=\sum\limits_{n\geq m}P_{n,m}(x)z^n$,
comparing the coefficients of $z^n$ on both side of Identity (2), we
obtain \begin{eqnarray*}P_{m,m}(x)&=&b_m(x),~and\\
P_{n,m}(x)&=&x\sum\limits_{i=0}^{m}\sum\limits_{j=i}^{n-m+i-1}b_{m-i}(x)b_{n-m+i-j-1}(x)P_{j,i}(x).
\end{eqnarray*}This complete the proof.
\end{proof}

\section{The refinement of $(n,m)$-flawed paths obtained by double ascent}

In this section, we will consider the refinement of $(n,m)$-flawed
paths obtained by double ascent and prove the value of $a_{n,m,k}$
have the Chung-Feller property on $m$. Define the generating
functions $A_m(x,z)=\sum\limits_{n\geq m}A_{n,m}(x)z^n$. When $m=0$,
$a_{n,0,k}=\frac{1}{k+1}{n-1\choose{k}}{n\choose{k}}$. It is well
know that $$A_0(x,z)=1+\frac{zA_0(x,z)}{1-xzA_0(x,z)},$$
equivalently,
$A_0(x,z)=\displaystyle{\frac{1+(x-1)z-\sqrt{(1+zx-z)^2-4xz}}{2xz}}.$

Define a generating function $A(x,y,z)=\sum\limits_{n\geq
0}\sum\limits_{m=0}^{n}\sum\limits_{k=1}^{n}a_{n,m,k}x^ky^mz^n.$
\begin{lem}\label{generatingdoubleascent}
$$A(x,y,z)=\displaystyle{\frac{A_0(x,z)A_0(x,yz)}{1-x[A_0(x,z)-1][A_0(x,yz)-1]}}.$$
\end{lem}
\begin{proof} Let the mapping $\phi$ be defined as that in
Introduction. An alternating Catalan path is a Dyck path which can
be decomposed into $RT$, where $R\neq\emptyset$ and $T\neq
\emptyset$, such that $\phi(R)$ and $T$ are both Catalan paths.

Now, Let $P\in \mathcal{D}$. We can uniquely decompose $P$ into
$PQ_1\ldots Q_mR$ such that $P$ and $\phi(R)$ are Catalan paths and
$Q_i$ is the alternating Catalan path for all $i$. Hence,
\begin{eqnarray*}A(x,y,z)&=&A_0(x,z)\left(\sum\limits_{m\geq 0}x[A_0(x,z)-1][A_0(x,yz)-1]\right)A_0(x,yz)\\
&=&\displaystyle{\frac{A_0(x,z)A_0(x,yz)}{1-x[A_0(x,z)-1][A_0(x,yz)-1]}}.\end{eqnarray*}
\end{proof}
\begin{thm}\label{chungfellerdoubleascent} Let $n$ be an integer with $n\geq 0$ and $0\leq k\leq
n-1$. Then
\begin{eqnarray*}a_{n,m,k}=\frac{1}{k+1}{n-1\choose{k}}{n\choose{k}}\end{eqnarray*}
for any $0\leq m\leq n$
\end{thm}

\begin{proof} First, we give an algebra proof of this theorem.
Since $xz[A_0(x,z)]^2=A_0(x,z)[1+xz-z]-1$, simple calculations tell
us
\begin{eqnarray*}
&&z\left\{1-x[A_0(x,z)-1][A_0(x,yz)-1]\right\}[yA_0(x,yz)-A_0(x,z)]\\
&=&z(y-1)A_0(x,z)A_0(x,yz).
\end{eqnarray*}
 By Lemma \ref{generatingdoubleascent}, we have
\begin{eqnarray*}
A(x,y,z)&=&\displaystyle{\frac{A_0(x,z)A_0(x,yz)}{1-x[A_0(x,z)-1][A_0(x,yz)-1]}}\\
&=&\displaystyle{\frac{yA_0(x,yz)-A_0(x,z)}{y-1}}\\
&=&\sum\limits_{n\geq 0}\sum\limits_{m=0}^nA_{n,0}(x)y^mz^n.
\end{eqnarray*}
This implies $A_{n,m}(x)=A_{n,0}(x)$ for any $0\leq m\leq n$.
Therefore,
$a_{n,m,k}=a_{n,0,k}=\frac{1}{k+1}{n-1\choose{k}}{n\choose{k}}$.

Now, we give a bijection proof of this theorem. Let $P$ be a path of
semi-length $n$ with $m$ flaws and $k$ double ascent, where $0\leq
m\leq n-1$. We say that a catalan path is prime if the path touches
$x$-axis exact twice. We can decompose $P$ into $SRUQDT$ such that

(1) $UQD$ is the right-most prime catalan path in $P$

(2) $\phi(R)$ is a catalan path, where $\phi$ is defined as that in
Introduction, and

(3) the final step of $S$ is $D$ on $x$-axis or $S=\emptyset$.\\
It is easy to see $\phi(T)$ is a catalan path. We define a path
$\varphi(P)$ as $$\varphi(P)=STDRUQ.$$ Clearly, the number of double
ascents in $\varphi(P)$ is equal to the number of double ascents in
$P$ and the number of flaws in $\varphi(P)$ is $m+1$.

To prove the mapping $\varphi$ is a bijection, we describe the
inverse $\varphi^{-1}$ of the mapping $\varphi$ as follows:\\
Let $P'$ be a path of semi-length $n$ with $m+1$ flaws and $k$
double ascent, where $0\leq m\leq n-1$. We can decompose $P'$ into
$STDRUQ$ such that

(1) $D$ and $U$ are the right-most steps leaving and returning to
the $x$-axis steps and under the $x$-axis in $P'$;

(2) $\phi(T)$ is a catalan path, where $\phi$ is defined as that in
Introduction, and

(3) the final step of $S$ is $D$ on $x$-axis or $S=\emptyset$.\\
Clearly, $Q$ and $\phi(DRU)$ are  both catalan paths. We define a
path $\varphi^{-1}(P')$ as  $\varphi^{-1}(P')=SRUQDT$.

\end{proof}
\begin{center}
\includegraphics[width=12cm]{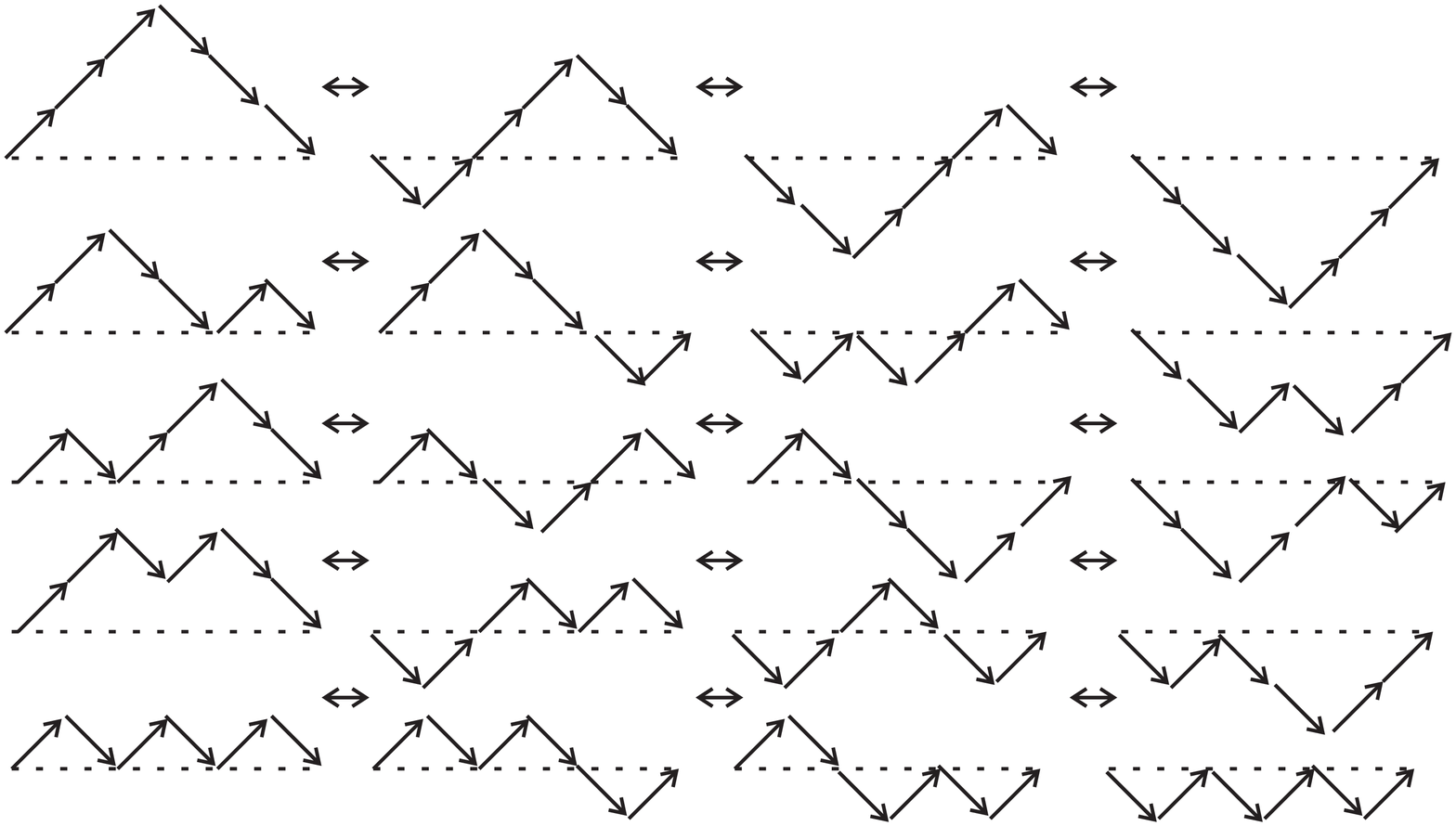}\\
Fig.2. $3$-Dyck path $P$ with $m$ flaws and $k$ double ascents and
$\varphi(P)$
\end{center}

\begin{cor}(Chung-Feller.) The number of $n$-Dyck path with
$m$-flaws is the Catalan number $c_n$ for any $0\leq m\leq n$.
\end{cor}

Recall that $d_{n,m,k}$ is the number of Dyck paths of semi-length
$n$ with $m$ flaws and $k$ double descents and
$d_{n,m,k}=a_{n,m,k}$.
\begin{cor}
 Let $n$ be an integer with $n\geq 0$ and $0\leq k\leq
n-1$. Then
\begin{eqnarray*}d_{n,m,k}=\frac{1}{k+1}{n-1\choose{k}}{n\choose{k}}\end{eqnarray*}
for any $0\leq m\leq n$.
\end{cor}

 \end{document}